\documentstyle{amsppt}
\pageno=1

\TagsOnRight
\magnification=\magstep1
\pageno1

\define\llapp#1{\hbox{\llap{$\dsize{#1}$}}}
\define\rlapp#1{\hbox{\rlap{$\dsize{#1}$}}}
\define\E{{\,{\bold E}\,}}
\redefine\P{\,{\bold P}\,}
\define\a{\alpha}
\define\be{\beta}
\define\g{\gamma}

\define\de{\delta}

\define\<{\bigl <}\define\8{\left<}\define\9{\right>}
\define\>{\bigr>}
\define\ovln#1{\,{\overline{\!#1}}}

\redefine\L{\Cal L}

\define\vk{\varkappa}

\define\bgm{\,\big|\4}

\define\si{\sigma}
\define\me{^{-1}}
\define\4{\kern1pt}
\define\wt{\widetilde}
\define\wh{\widehat}
\define\2{\kern.5pt}
\define\6{\phantom0}
\define\bgr#1{\4\bigr#1}
\define\bgl#1{\bigl#1\4}
\define\3{\kern-10pt}
\define\half{{}^1\kern-.5pt\!\big/\!\2{}_2}
\define\sfrac#1#2{{}^{#1}\kern-.5pt\!\big/\!\kern.5pt{}_{#2}}
\define\pfrac#1#2{{#1}\big/{#2}}
\define\ffrac#1#2{\raise.5pt\hbox{\eightpoint$\4\dfrac{\,#1\,}{\,#2\,}\4$}}
\define\Bgr#1{{\eightpoint\4\raise.5pt\hbox{{$\Bigr#1$}}}}
\define\Bgl#1{{\eightpoint\raise.5pt\hbox{{$\Bigl#1$}}\4}}
\redefine\le{\leqslant}
\redefine\ge{\geqslant}
\def\={\overset{ \text{def} }\to =}
\define\na{\,\, {\raise.4pt\hbox{$\shortmid$}}{\hskip-2.0pt\to}\, \, }
\define\BR{\bold R}

\hoffset .25 true cm
\hsize=15.5 true cm
\vsize=21 true cm

\document
\topmatter
\title A multiplicative inequality for  concentration functions
of $n$@-fold convolutions
\endtitle
\rightheadtext{ Estimates for   concentration functions}
\leftheadtext{  F. G\"otze and A. Yu. Zaitsev}  
\author F. G\"otze$^1$\quad A. Yu. Zaitsev$^{1,2}$
\endauthor
\thanks
{$^1$Research supported  by the SFB 343 in Bielefeld.
\newline\indent
$^2$Research supported by Russian Foundation of Basic Research
(RFBR) Grant 05-01-00911,  by RFBR-DFG Grant 04-01-04000,
and  by the grant NSh 4222.2006.1}
\endthanks
\affil 
University of Bielefeld$^1$\\
St.~Petersburg Branch of Steklov Mathematical Institute$^2$
 \endaffil
\address  
Friedrich G\"otze \smallskip
Fakult\"at f\"ur Mathematik \smallskip
Universit\"at Bielefeld \smallskip
Postfach 100131 \smallskip
33501 Bielefeld 1 \smallskip
Germany  \smallskip
\endaddress
\email
goetze\@mathematik.uni-bielefeld.de
\endemail
\address Andrei Yu.  Zaitsev\smallskip
St.~Petersburg Branch of Steklov Mathematical Institute
\smallskip
 Fontanka 27\smallskip St.~Petersburg 191011\smallskip Russia
\smallskip
\endaddress
\email
zaitsev\@pdmi.ras.ru
\endemail
\date   December 1999  
\enddate
\abstract
We estimate  
the concentration functions of $n$@-fold convolutions
of one-di\-men\-sional probability measures.
 The main result is a supplement
to the results of G\"otze and Zaitsev (1998).
We show that the estimation of concentration functions 
at arguments of bounded size can be reduced  to
the estimation of these functions at arguments of size
$O(\sqrt n)$ which is easier.
\endabstract
\endtopmatter

\head{\bf1. Introduction}\endhead
Let us first introduce some notation. Let
$\goth F$ denote the set of probability 
distributions defined on the
Borel $\si$-field of subsets of the real line~$\BR $, 
\,$\L(\xi)\in\goth F$ \,
the distribution of a random variable~$\xi$, 
\,and \,$\bold I\{A\}$ \,the
indicator function of an event~\,$A$. 
For $F\in\goth F$  the concentration function 
is defined by 
\;${Q(F,\,b)=\sup_x F\bgl\{[x,\,x+b]\bgr\}}$, \, $b\ge0$.
For $F,H\in\goth F$ we
denote the corresponding  distribution functions by~$F(x)$,~$H(x)$ \,
and the characteristic functions by ~$\wh F(t)$,~$\wh H(t)$. \,
Let \,$E_a\in\goth F$ \,be the
distribution concentrated at a point $a\in \BR $, $E=E_0$.
Products and powers of measures 
will be understood in the convolution
sense: \; ${FH=F*H}$, \,$H^m=H^{m*}$, \,$H^0=E$. 
For \,$F=\L(\xi)\in\frak F$ \, we shall use the notation
 \,$\ovln F=\L(-\xi)$ \,and  \,$\wt F=F\4\ovln F$. \,
The distribution~ \,$\wt F$ \, is  called symmetrized.
By~\, $c(\4\cdot\4)$ \,
we shall  denote   positive quantities
depending on the arguments only. The
 symbol~\,$c$ \,will be used for 
absolute positive constants. Note that \,$c(\4\cdot\4),\,c$ \,
can be different  in different
(or even in the same) formulas. 
The ends of proofs will be marked by ~\,$\square$.

We estimate  
the concentration functions of $n$@-fold convolutions
of one-di\-men\-sional probability measures.
The concentration functions was introduced 
and preliminarily estimated by
L\'evy (1937), see as well Doeblin (1939).
The bounds for these convolutions
were obtained by  Kolmogorov (1956,~1958), Rogozin~(1961a,b), 
Le Cam (1965),
 Esseen~(1966,~1968), Kesten (1969, 1972), Morozova (1977),
 Postnikova and Yudin (1978),
 Arak~(1981), Miroshnikov and Rogozin (1980, 1982),
Hall ~(1983), Griffin, Jain and Pruitt (1984),
Arak and Zaitsev (1988), Nagaev and  Khodzhabagyan (1996)
and others.

The aim of this paper is to provide a supplement
to a previous paper G\"otze and Zaitsev (1998)
which is abbreviated in the following as GZ. 
We generalize Theorem~ 2 of that paper proving
the following result.

\proclaim{Theorem~1}
Assume that a distribution
\,$F\in\goth F$ \,
is represented in the form
$$
F=(1-p)\4U+p\4V,
\qquad U,V\in\goth F, \quad 0<p<1.
\tag1.1
$$
Let \,$X$ \,be a random variable with \;$\L(X)=U$.
\;Suppose that
$$
0<\si^2=\E X^2<\infty,
\qquad
\E X=0,
\tag1.2
$$
and
$$
 b\ge\ffrac{\vk_n}{\si^2}, \qquad
\text{where}\quad
\vk_n=\E
 X^2\4\min\bgl\{|X|,\,\si\sqrt n\bgr\} .
\tag1.3
$$
Let \,$r, s$ \,be integers, \,$0\le r<s\le n$. \,Then,
for any distribution \,$H\in\goth F$,
$$
\multline
 Q(H\4F^n,\,b)\le
\E\ffrac{c\4b}{\si\sqrt {n-\mu}}\, 
Q\big(H\4V^\mu,\,\si\sqrt {n-\mu}\big)\,
\bold I\bgl\{r\le\mu<s\bgr\}\\
+\E Q\big(H\4U^{n-\mu}\4V^\mu,\,b\big)\,
\bold I\bgl\{\mu\ge s\bgr\}
+\min\Bigl\{\,1,\,\ffrac{c\4b}{\si\sqrt {n-r}}
\,Q\big(H,\,\si\sqrt {n}\big)\Bigr\}\P\bgl\{\mu<r\bgr\},
\endmultline
\tag1.4
$$
where \,$\mu$  \,is a random variable
having binomial distribution with parameters \,$n$ \,and ~\,$p$. 
\,Moreover,
$$
\multline
 Q(H\4F^n,\,b)\le
\ffrac{c\4b}{\si\sqrt{n\,(1-p)}}\, 
Q\big(H\4V^r,\,\si\sqrt {n}\big)\\
 +
\min\Bigl\{\,1,\,\ffrac{c\4b}{\si\sqrt {n-r}}
\,Q\big(H,\,\si\sqrt {n}\big)\Bigr\}\P\bgl\{\mu<r\bgr\}.
\endmultline
\tag1.5
$$
\endproclaim

Theorem~2 of GZ (1998)
is a particular case of our Theorem ~1 appearing
when  \,$H=E$. 
 The moment restrictions
 in Theorem~1 are imposed
 on the distribution~$U$ only.
The distributions $H$ and $V$  are  arbitrary
and therefore the initial distribution
~\,$F$ \,may have arbitrarily bad moment properties.

Taking in (1.5) \,$r=0$, \,we obtain

\proclaim{Corollary 1}Under the conditions of Theorem $1$
$$
 Q(H\4F^n,\,b)\le
\ffrac{c\4b}{\si\sqrt{n\,(1-p)}}\, 
Q\big(H,\,\si\sqrt {n}\big).
\tag1.6
$$
\endproclaim

An  important particular case of Theorem 1  and  Corollary 1 
appears when \,$n$ \,is changed by integer \,${n\4\a}$ \,and
 \,$H=F^{n(1-\a)}$ \,where \,$0<\a<1$. \, Then 
the bound (1.6) can be rewritten in the form
$$
 Q(F^n,\,b)
\le
\ffrac{c\4b}{\si\sqrt{n\4\a\4(1-p)}} \,
Q\big(F^{n(1-\a)},\,\si\sqrt {n\4\a}\big). 
\tag1.7
$$
Theorem 1  and  Corollary 1 
provide in a sense
"multiplicative inequalities" for concentration functions of
convolutions. The inequality (1.7) 
reduces the problem of estimation
of~\,$Q(F^n,\,b)$ \,to that of 
of~\,$Q\big(F^{n(1-\a)},\,\si\sqrt {n\4\a}\big)$. \,
The trivial inequality
$$
Q\big(F^{n(1-\a)},\,\si\sqrt {n\4\a}\big) 
\le  \Bigl(1+\ffrac{\si\sqrt{n\4\a}}{b}\Bigr)\, 
Q(F^{n(1-\a)},\,b)
\tag1.8
$$
(see~(2.4)) together with~(1.7) may return the problem
to the initial stage (in particular, when
\,$Q(F^n,\,b)$ \,
decreases with some negative power of~\,$n$).
However, the inequality~(1.8) may lead
to a loss of precision. For instance, if
\,$F$ \, is concentrated on a 
one-dimensional lattice with step size \,$2\4\si\sqrt{n\4\a}$, \,
 we have
$$
Q\big(F^{n(1-\a)},\,\si\sqrt {n\4\a}\big) 
= Q(F^{n(1-\a)},\,b), \qquad\text{for} \quad
0\le b\le\si\sqrt{n\4\a}.
$$
It is clear that it is much easier to estimate the concentration
function for the large value \,$\si\sqrt{n\4\a}$ of the argument than for 
some fixed ~\,$b$. \, 

Applying Theorem~1 and Corollary 1, 
one can use known classical bounds for concentration functions.
We mention in this connection
the papers by Rogozin~(1961a), 
 Esseen~(1968), \;Kesten (1969),
\;Arak~(1981), \; Miroshnikov \;and \;Rogozin (1980),
Hall ~(1983), Griffin, Jain and Pruitt (1984).
For example,  the Esseen bound (1968), 
implies the following
result.

\proclaim{Corollary~2}Let the conditions of Theorem~$1$
be satisfied. Then, for any \, $\de>0$,
$$
 Q(F^n,\,b)\le
\ffrac{c\4b\,(\de+\si)}{\de\4\si\4n\4
\sqrt{(1-p)\,D\bgl(\wt F,\,\de\sqrt{n}\bgr)}}.
\tag1.9
$$
 where
$$
 D(F,\,b)=
\int_{-\infty}^\infty \min\bgl\{x^2\4b^{-2},\,1\bgr\}\,F\{dx\},
\qquad F\in\frak F,\quad b>0.
\tag1.10
$$
\endproclaim
 
Esseen~(1968)  (see Petrov~(1976), inequality~(2.7)
of Chapter~III)
proved that
$$
 Q(F^n,\,b)\le
\ffrac{c}{\sqrt{n\4D\bgl(\wt F,\,b\bgr)}},
\tag1.11
$$
Corollary ~2 implies the following result.

\proclaim{Corollary~3}
Let \,$F\in\frak F$, \,$b\ge0$, \,$\de>0$. \,Then
$$
 Q(F^n,\,b)\le
\ffrac{c(F,b,\de)}{n\4\sqrt{D\bgl(\wt F,\,\de\sqrt{n}\bgr)}}.
$$
\endproclaim

For the proof it suffices to note that after a shift
every non-degenerate distribution can be 
represented in the form (1.1) with \,$p=\half$ \, 
and a non-degenerated $U$
having bounded support. This yields the result for
\,$b\ge \vk/\si^2$, \,where \,$\vk=\E|X|^3$. \,
For  \,$b< \vk/\si^2$ \,we can apply the result
for  \,$b= \vk/\si^2$, \,using the monotonicity of
concentration functions. 
Corollary~3 is Theorem~1 from the paper of
GZ (1998) which is a sharpening of a result of
Esseen~(1968) who  showed  that  \;$Q(F^n,\,b)=o(n^{-1/2})$ \,
as  \,$n\to\infty$ \,
iff the distribution~\,$F$ \,has an infinite second moment
(see as well Morozova~(1977)).
However the proof of Corollary~3
in the present paper is somewhat easier
and we obtain here a more explicit form of~\,$c(F,b,\de)$. \,
 For the connection of Corollary~3
with previous results about concentration functions see GZ (1998).

Note that in the proof of  Corollary~3
we use the inequality~(1.11). Nevertheless,
Corollary~3 can be considered as an improvement
of the inequality~(1.11). Indeed,  the latter 
can be rewritten in the form
$$
 Q(F^n,\,b)\le
c\,\biggl(n\4
\int\min\bgl\{x^2\4b^{-2},\,1\bgr\}\,
\wt F\{dx\}\biggr)^{\!-1/2}
\tag1.12
$$ 
Comparing (1.12) with the inequality
$$
 Q(F^n,\,b)\le
c(F,b)\,\biggl(n\4
\int\min\bgl\{x^2\4b^{-2},\, n\bgr\}\,
\wt F\{dx\}\biggr)^{\!-1/2}.
\tag1.13
$$ 
which follows from Corollary~3 with \,$\de=b$, \,we see that for
any distributions with  infinite variance
the inequality (1.13) is sharper with respect to the order in ~\,$n$ \,
 than~ (1.12). Note however that it is impossible
to change $c(F,b)$ in (1.13) by some absolute constant~\,$c$.
The corresponding example is given by the distribution
\,$F=F_n=\ffrac12E_{-n}+\ffrac12E_{n}$. \,
It is well known that \,$Q(F_n^n,1)$ \,
behaves as \,$O(n^{-1/2})$ \,when \,$n\to\infty$. \,
On the other hand,
$$
\biggl(n\4
\int\min\bgl\{x^2,\, n\bgr\}\,
\wt F_n\{dx\}\biggr)^{\!-1/2}\le \ffrac cn.
$$ 
This implies that \,$ c(F_n,1)\ge c\4\sqrt n$ \,in (1.13).

Esseen (1968) proved that, for any \,$F\in\goth F$ and \,$b>0$,
$$
Q(F,b)\le c\4b\,\int\limits_{|t|\le b\me}
\bgl|\wh F(t)\bgr|\,dt.
\tag1.14
$$
Applying (1.14) to the distribution \,$F^n$, \,
we obtain
$$
Q(F^n,b)\le c\4b\,\int\limits_{|t|\le b\me}
\bgl|\wh F(t)\bgr|^n\,dt.
\tag1.15
$$
On the other hand, using (1.7) 
and then (1.14), we see that
$$
Q(F^n,b)\le \ffrac{c\4b}{\sqrt{1-p}}\,
\int\limits_{|t| \si\sqrt{n\a}\le1}
\bgl|\wh F(t)\bgr|^{n(1-\a)}\,dt,\nopagebreak
\tag1.16
$$
if \,$n\4\a$ \,is integer and \,$b\ge\E|X|^3/\si^2$. 
\,A comparison of (1.15) with  (1.16) demonstrates
an advantage of using the inequality~(1.7).

Our proofs  are based on  non-uniform
estimates in the Central Limit Theorem (CLT)
and on elementary properties of concentration functions
(see the proof of Lemma~1 in GZ (1998) and 
Zaitsev (1987, ~1992)).
In this respect our proofs
 differ  from most of the previous
papers, where Esseen's (1968) method of 
characteristic functions had been extensively used.
One should note however that the CLT approach 
was applied in the seminal paper
of Kolmogorov (1958). He used the uniform
Berry--Esseen bound in the CLT since
 non-uniform
ones  were not known at that time.

\head{\bf2. Proofs}\endhead

\proclaim{Lemma 1 \rm(GZ (1998, Lemma 3))} 
Let \,$\xi,\xi_1,\xi_2,\dots,\xi_n$
 \,be i.i.d\. random variables,
\,$\E\xi=0$, \;
$ W=\L(\xi_1+\dots+\xi_n)$, 
$$
B^2=n\E\xi^2>0, \qquad
\be=n\E \xi^2\4\min\bgl\{|\xi|,\,B\bgr\}<\infty .
\tag2.1
$$
and let
$
b\ge\ffrac{\be}{B^2}
$. \,
 Then, for any  \,$G\in\goth F$, \,
we have
$$
Q(W G,\,b)\le\ffrac{c\4b}{B}\4Q(G,B).
\tag2.2
$$
\endproclaim

Note that the inequality (2.2) is a particular case
of the inequality (1.6) appearing when \,$p=0$.

\smallbreak

We need the following
well-known simple properties of 
concentration functions, which are valid for any
$F,H\in\goth F$ and
$\g,\g_1,\g_2>0$:
$$\spreadlines{1\jot}
\align
Q(FH,\g) &\le\min\bgl\{Q(F,\g),\,Q(H,\g)\bgr\};
\tag2.3 \\
Q(F,\g_1) &\le\bgl (1+\lceil\g_1/\g_2\rceil\bgr )\,Q(F,\g_2),
\tag2.4
\endalign
$$
where
$\lceil\4\cdot\4\rceil$  is the integer part of a number
(see, e.g.,  Hengartner and  Theodorescu~(1973)).

\pagebreak

\demo{Proof of Theorem $1$}
It is known that 
$H\4F^n$ can be written in the form
$$
H\4F^n=\L (S_n ),
\qquad\text{where}\quad
S_n=\zeta+\sum_{i=1}^n(1-\mu_i)\4\xi_i+\mu_i\4\eta_i,
$$
and \,$\zeta$, \,$\xi_i$, $\eta_i$, $\mu_i$ \,
are jointly independent random
variables with 
$$
\L(\zeta)=H,\qquad
\L(\xi_i)=U,\qquad \L(\eta_i)=V, \qquad
 \L(\mu_i)=(1-p)\4E+p\4E_1.
$$
Define
\;$\mu=\sum\limits_{i=1}^n\mu_i$. \;
Obviously, \,$\mu$  \,has binomial 
distribution with parameters~\,$n$ \,and~\,$p$. \,
Given a fixed value of~\,$\mu$, \,the random
variable~\,$S_n$ \,has  conditional distribution
~\,$H\4U^{n-\mu}\4V^\mu$. \,
Hence, for any \,$x\in\bold R$ \,we have
$$
\split
\P\Bgl\{S_n\in\bgl[x,\,x+b\bgr]\Bgr\}&=
\E\,\bold I\bgl\{S_n\in\bgl[x,\,x+b\bgr]\bgr\}\\
&= \E\E \Bgl\{\bold I\bgl\{S_n\in\bgl[x,\,x+b\bgr]\bgr\}
 \bgm  \mu\Bgr\}\\
&\le\E Q\big(H\4U^{n-\mu}\4V^\mu,\,b\big).
\endsplit
$$
Therefore,
$$
Q(H\4F^n, b)\le\E Q\big(H\4U^{n-\mu}\4V^\mu,\,b\big).
\tag2.5
$$
In view of \,$\vk_{n-\mu}\le\vk_n$
\,and 
applying inequality~(2.2) of Lemma~1 with
\;$W=U^{n-\mu}$,
\;$G=H\4V^\mu$,
\;$B=\si\sqrt {n-\mu}$, \;
 we see that
$$
\multline
\E Q\big(H\4U^{n-\mu}\4V^\mu,\,b\big)\le
\E\ffrac{c\4b}{\si\sqrt {n-\mu}}\, 
Q\big(H\4V^\mu,\,\si\sqrt {n-\mu}\big)\,
\bold I\bgl\{r\le\mu<s\bgr\}\\
+\E Q\big(H\4U^{n-\mu}\4V^\mu,\,b\big)\,
\bold I\bgl\{\mu\ge s\bgr\}
+\E Q\big(H\4U^{n-\mu}\4V^\mu,\,b\big)\,
\bold I\bgl\{\mu<r\bgr\}
\endmultline
\tag2.6\nopagebreak
$$
and
$$
\E Q\big(H\4U^{n-\mu}\4V^\mu,\,b\big)\,\bold I\bgl\{\mu<r\bgr\}\le
\min\Bigl\{\,1,\,\ffrac{c\4b}{\si\sqrt {n-r}}\,
Q\big(H,\,\si\sqrt {n}\big)\Bigr\}
\P\bgl\{\mu<r\bgr\}.
\tag2.7\nopagebreak
$$
The inequality (1.4) now follows from (2.5)@--(2.7).

Using the relations \,$\P\bgl\{\mu=n\bgr\}=p^n$,
\, (1.4) with \,$s=n$ \,and ~(2.3), we obtain
$$
\multline
Q(H\4F^n,\,b)\le
\ffrac{c\4b}{\si}\, 
Q\big(H\4V^r,\,\si\sqrt {n}\big)
\E\ffrac{1}{\sqrt {n-\mu}}\, 
\bold I\bgl\{r\le\mu<n\bgr\}\\
+p^n\,Q\big(H\4V^n,\,b\big)+
\min\Bigl\{\,1,\,\ffrac{c\4b}{\si\sqrt {n-r}}
\,Q\big(H,\,\si\sqrt {n}\big)\Bigr\}
\P\bgl\{\mu<r\bgr\}.
\endmultline
\tag2.8
$$
Applying the H\"older inequality, we derive 
$$
\split
\llapp{\E\ffrac{1}{\sqrt {n-\mu}}}\, 
\bold I\bgl\{r\le\mu<n\bgr\}
&\le c\,\E\ffrac{1}{\sqrt {n-\mu+1}}\\
&\le c\,\Bigl(\E\ffrac{1} {n-\mu+1}\Bigr)^{\!1/2}
\le \rlapp{\ffrac{c} {\sqrt{n\,(1-p)}}.}
\endsplit
\tag2.9
$$
Using~(1.2) and~(1.3), it is easy to see that
$$
\split
\llapp{\vk_n}&\ge\vk_1=
\E|X|^3\,\bold I\bgl\{|X|\le\si\bgr\}
+\si\E|X|^2\,\bold I\bgl\{|X|>\si\bgr\}
\\
&\ge\Bgl(\E|X|^2\,\bold I\bgl\{|X|\le\si\bgr\}\Bgr)^{\!3/2}
+\Bgl(\E|X|^2\,\bold I\bgl\{|X|>\si\bgr\}\Bgr)^{\!3/2}
\ge \rlapp{c\4\si^3.}
\endsplit
\tag2.10
$$
The relations (1.3) and (2.10) together imply
$$
b\ge\vk_n\4\si^{-2} \ge c\4\si.
\tag2.11
$$
According to~(2.3), (2.4) and (2.11), we have
$$
\split
 Q\big(H\4V^n,\,b\big)
\le Q\big(H\4V^r,\,b\big)
&\le\Bigl(\ffrac{b}{\si\sqrt {n}}+1\Bigr)
\,Q\big(H\4V^r,\,\si\sqrt {n}\big)\\
&\le\ffrac{c\4b}{\si}
\,Q\big(H\4V^r,\,\si\sqrt {n}\big).
\endsplit
\tag2.12
$$
Moreover,
$$
p^n=\big(1-(1-p)\big)^n\le
e^{-n\2(1-p)}\le \ffrac{c} {\sqrt{n\,(1-p)}}.
\tag2.13
$$
The inequality (1.5) can be easily
derived  from~(2.8), (2.9), (2.12) and (2.13).
$\square$\enddemo\smallbreak

\demo{Proof of Corollary~$2$}
In view of~(2.3), (2.11),
 we may assume without loss of generality that \,$n$ 
is even and \,$n\ge2$.
\,Using  (1.7) with \,$\a=\half$, \,(1.11), (2.3)  and
(2.4), we obtain
$$
\split
 Q(F^n,\,b)&\le
\ffrac{c\4b}{\si\sqrt{n\4(1-p)}}\, 
Q\bgl(F^{n/2},\,\si\sqrt {n}\bgr) \\
&\le
\ffrac{c\4b}{\si\sqrt{n\4(1-p)}}\,\bgl(1+\si\4\de\me\bgr)\, 
Q\bgl(F^{n/2},\,\de\sqrt {n}\bgr)\\
&\le
\ffrac{c\4b}{\si\sqrt{n\4(1-p)}}\,
\ffrac{1+\si\4\de\me}{\sqrt{n\,D\bgl(\wt F,\,\de\sqrt{n}\bgr)}}\\
&=
\ffrac{c\4b\,(\de+\si)}
{\de\4\si\4n\4\sqrt{(1-p)\,D\bgl(\wt F,\,\de\sqrt{n}\bgr)}}.
\qquad_{\dsize\square}
\endsplit
$$
\enddemo

\newpage

\enddocument